# Rectangularization of Gaussian process regression for optimization of hyperparameters

Sergei Manzhos[1] and Manabu Ihara

School of Materials and Chemical Technology, Tokyo Institute of Technology, Ookayama 2-12-1, Meguro-ku, Tokyo 152-8552 Japan. E-mail: mihara@chemeng.titech.ac.jp (M.I.), manzhos.s.aa@m.titech.ac.jp (S.M.)

**Abstract**

Gaussian process regression (GPR) is a powerful machine learning method which has recently enjoyed wider use, in particular in physical sciences. In its original formulation, GPR uses a square matrix of covariances among training data and can be viewed as linear regression problem with equal numbers of training data and basis functions. When data are sparse, avoidance of overfitting and optimization of hyperparameters of GPR are difficult, in particular in high-dimensional spaces where the data sparsity issue cannot practically be resolved by adding more data. Optimal choice of hyperparameters, however, determines success or failure of the application of the GPR method. We show that parameter optimization is facilitated by rectangularization of the defining equation of GPR. On the example of a 15-dimensional molecular potential energy surface we demonstrate that this approach allows effective hyperparameter tuning even with very sparse data.

## 1 Introduction

The use of the Gaussian process regression (GPR) (Bishop, 2006; Rasmussen and Williams, 2006) approach is a powerful machine learning tool. GPR is easy to use, in particular, in high-dimensional spaces: being a non-parametric method (only a small number hyperparameters need to be selected), the increase in space dimensionality does not lead to a drastic increase in the number of (non-linear) parameters. This is in contrast to, for example, neural networks (NN) (Emmert-Streib et al., 2020; Montavon et al., 2012), where the number of nonlinear parameters rapidly proliferates when the numbers of neurons, layers, and space dimensionality increase. It has been argued that GPR's

---

[1] Author to whom correspondence should be addressed. E-mail: E-mail: manzhos.s.aa@m.titech.ac.jp , Tel & Fax : +81-3-5734-3918

expressive power is equivalent to that of an infinite-width neural network (NN) (Neal, 1995; Williams, 1996). This increases the danger of overfitting as well as the calculation cost. Recent appearance of comparisons between NN and GPR on the same data and in applications requiring very high accuracy such as fitting of spectroscopically accurate potential energy surfaces (PES) (Kulik and et al., n.d.; Majumder et al., 2015; Manzhos et al., 2015b, 2006; Manzhos and Carrington, 2021) and kinetic energy functionals (Manzhos and Golub, 2020) highlighted GPR advantages, in particular, in obtaining highly accurate approximations of multidimensional functions from few data (Kamath et al., 2018).

One problem area, as will also be highlighted below, is overfitting and the optimal choice of GPR hyperparameters (of the GPR kernel) that allows to avoid it, which is especially difficult when data are scarce. While in high-dimensional spaces GPR is especially attractive, as growth of the space dimensionality need not lead to increased computational cost of the key matrix equations 1 and 2 below, high-dimensional applications need to deal with the issue of hyperparameter optimization under low density of training data. We note that data are always sparse in sufficiently high-dimensional spaces (Donoho, 2000). The problem of GPR hyperparameter optimization has been addressed with various methods; examples are random search related methods (Bergstra and Bengio, 2012), various versions of Bayesian inference (Brochu et al., 2010; Snoek et al., 2012) (with the commonly used maximum likelihood estimator, also used in this work, belonging to this type of methods (Myung, 2003)), simulated annealing (Fischetti and Stringher, 2019), genetic algorithms,(Alibrahim and Ludwig, 2021) so-called bandit-based methods (Li et al., 2018), and combinations thereof (Falkner et al., 2018). However, the scarcity (low density) of available training (and test) data may not permit optimizing the kernel parameters for the best *global* (i.e. in all relevant parts of the descriptors space) quality of function representation (i.e. to avoid overfitting). The global quality of representation can be considered from the point of view of the completeness of the basis set when one views GPR as linear regression with a basis formed by the kernel functions (see below).

To briefly summarize GPR (Rasmussen and Williams, 2006), we consider the problem of reconstructing a continuous function $f(x), x \in R^D$ from a finite number of samples $f^{(j)} = f(x^{(j)})$ at points $x^{(j)}, j = 1, \ldots, M$. The expectation values $f(x)$ and variances $\Delta f(x)$ of the functions values at any point in space $x$ are computed as (Rasmussen and Williams, 2006)

$$f(x) = K^* K^{-1} f$$

(1)

$$\Delta f(x) = K^{**} - K^* K^{-1} K^{*T}$$

(2)

where $f$ is a vector of $f^{(j)}$, i.e. the $M$ known values of the target function, the $M \times M$ matrix $\boldsymbol{K}$ is computed from pairwise covariances among the data:

$$\boldsymbol{K} = \begin{pmatrix} k(\boldsymbol{x}^{(1)}, \boldsymbol{x}^{(1)}) + \delta & k(\boldsymbol{x}^{(1)}, \boldsymbol{x}^{(2)}) & \cdots & k(\boldsymbol{x}^{(1)}, \boldsymbol{x}^{(M)}) \\ k(\boldsymbol{x}^{(2)}, \boldsymbol{x}^{(1)}) & k(\boldsymbol{x}^{(2)}, \boldsymbol{x}^{(2)}) + \delta & & k(\boldsymbol{x}^{(2)}, \boldsymbol{x}^{(M)}) \\ \vdots & & \ddots & \vdots \\ k(\boldsymbol{x}^{(M)}, \boldsymbol{x}^{(1)}) & k(\boldsymbol{x}^{(M)}, \boldsymbol{x}^{(2)}) & \cdots & k(\boldsymbol{x}^{(M)}, \boldsymbol{x}^{(M)}) + \delta \end{pmatrix}$$

(3)

the row vector $\boldsymbol{K}^*$ of length $M$ is

$$\boldsymbol{K}^* = \begin{pmatrix} k(\boldsymbol{x}, \boldsymbol{x}^{(1)}) & k(\boldsymbol{x}, \boldsymbol{x}^{(2)}) & \ldots & k(\boldsymbol{x}, \boldsymbol{x}^{(M)}) \end{pmatrix},$$

(4)

and $K^{**} = k(\boldsymbol{x}, \boldsymbol{x})$. The covariance function $k(\boldsymbol{x}^{(1)}, \boldsymbol{x}^{(2)}|\lambda)$ is the kernel of GPR that depends on hyperparameters $\lambda$ (which we generally omit in the formulas for brevity). The optional $\delta$ on the diagonal is a regularization (hyper)parameter; it is introduced to achieve numeric stability of the inversion of $\boldsymbol{K}$ and to improve the generalization of the model. Commonly used kernels belong to the Matern family (Genton, 2001):

$$k(\boldsymbol{x}, \boldsymbol{x}') = \sigma^2 \frac{2^{1-\nu}}{\Gamma(\nu)} \left( \sqrt{2\nu} \frac{|\boldsymbol{x} - \boldsymbol{x}'|}{l} \right)^\nu K_\nu \left( \sqrt{2\nu} \frac{|\boldsymbol{x} - \boldsymbol{x}'|}{l} \right)$$

(5)

where $\Gamma$ is the gamma function, and $K_\nu$ is the modified Bessel function of the second kind. At different values of $\nu$, this function becomes a squared exponential ($\nu \to \infty$), a simple exponential ($\nu = 1/2$) and other types of kernels (such as Matern3/2 and Matern5/2 for $\nu = 3/2$ and $5/2$, respectively). The value of $\nu$ is often preset, and the length scale $l$ and prefactor $\sigma^2$ are hyperparameters (i.e. $\lambda = (l, \sigma^2)$) that can be optimized (as only the ratio of the values of $\delta$ and $\sigma^2$ is important, in principle, only one of them needs to be optimized). In this work we will use squared exponential kernels with $\sigma^2 = 1$:

$$k(\boldsymbol{x}, \boldsymbol{x}') = \exp\left( -\frac{|\boldsymbol{x} - \boldsymbol{x}'|^2}{2\exp(l)^2} \right)$$

(6)

(note that in this case and if the values of $f(\boldsymbol{x})$ are not normalized to unit standard deviation, the right hand side of Eq. 2 needs to be multiplied by the variance of the target values $f$).

Overfitting is, however, a problem with GPR as well. An appropriate choice of hyperparameters is necessary to avoid it. We have recently shown (Manzhos and Ihara, 2022) that the commonly used maximum likelihood estimation (MLE) (Myung, 2003) that maximizes the log likelihood function,

$$\max\left(\frac{1}{2}ln|\mathbf{K}| - \frac{1}{2}\mathbf{t}\mathbf{K}^{-1}\mathbf{t} - \frac{N}{2}\ln(2\pi)\right)$$

(7)

can fail to find appropriate hyperparameters when data are sparse. While in principle it is tempting to think that this issue could be resolved by adding more data, in practice, this is a dead end proposition: in high-dimensional spaces, not only does adding more data not significantly change data density (in terms of the number data of data points per dimension of a corresponding direct product grid), the cost of wielding Eqs. 1-2 grows rapidly with $M$. We have shown in (Manzhos and Ihara, 2022) that even when using an additive model of $f(\mathbf{x})$, $f(\mathbf{x}) \approx \sum_{i=1}^{D} f_i(x_i)$ (Duvenaud et al., 2011; Manzhos et al., 2022), where the component functions $f_i(x_i)$ are well-determined (no overfitting) with few data (in the example of (Manzhos and Ihara, 2022), $f_i(x_i)$ were well-determined with only 500 data in a 15-dimensional space – precisely because they are low-dimensional as opposed to a 15-dimensional $f(\mathbf{x})$), it took about 10,000 training data for the MLE method to be able to find good hyperparameters under fixed $\delta$. Even with 10,000 training data, MLE did not find acceptable hyperparameters for a general non-additive GPR model or an additive model with simultaneous optimization of $l$ and $\delta$. This is in spite of the fact that hyperparameters do exist (and can be found manually) that result in an accurate GPR model without noticeable overfitting from 5,000 data (less than 1.8 data per dimension) (Manzhos and Ihara, 2022).

In this article, we explore an alternative way to determine optimal hyperparameters based on the view of GPR as a regularized linear model with a basis set made of the covariance functions. We view the problem of hyperparameter optimization as a basis completeness problem. In (Manzhos and Ihara, 2022), we proposed using RS-HDMR (random sampling high dimensional model representation) (Li et al., 2006; Rabitz and Aliş, 1999) to build, without the danger of overfitting, reference functions computable in all space which could be used to tune the hyperparameters (of the basis) by computing more data. That method has disadvantages that the basis completeness is optimized for a different target function than the function of interest (albeit close to it) and that more data need to be generated, increasing the cost of wielding the matrix equations of GPR. Here, we explore an alternative route of optimizing the basis completeness by modifying the governing equation of GPR, namely, by rectangularizing it.

## 2 Methods

We begin by recognizing that the equivalency between GPR and regularized linear regression (Bishop, 2006). Eq. 1 has the form of a basis expansion,

$$f(x) = \sum_{n=1}^{M} b_n(x) c_n$$

(8)

with basis functions $b_n(x) = k(x, x^{(n)})$ and linear coefficients $c_n$. The quality of this approximation is fully determined by the quality of the basis set $\{k(x, x^{(n)}|\lambda)\}$, i.e. the extent of its completeness. Ultimately it is the quality of the basis (in relevant parts of space) that hyperparameter optimization should improve. In the case of a squared exponential kernel used here, it is the quality of a Gaussian-type basis set with basis functions located at each of the $M$ training points.

Eq. 8 can be written in matrix form

$$f' = Bc$$

(9)

where the lengths of the vector $f'$ corresponds to the number $N$ of considered points $x'$. This could be all or a subset of points where the knowledge of $f(x)$ is required in an application (e.g. in the case of learning a molecular potential energy surface – which is the example of GPR application considered in this work – this could be points on a quadrature grid in a quantum dynamics calculation (Beck et al., 2000; Bowman et al., 2008)) or it could be a set of test points. $B$ is then a rectangular matrix of size $N \times M$ with elements $B_{mn} = k(x'^{(n)}, x^{(m)})$. If $\{x'\} = \{x\}$ (the same $M$ known points) then $B = K$ is a square matrix, $f = f'$ and $c$ is obtained as $c = K^{-1}f$ i.e. one obtains standard GPR (Bishop, 2006). The specific nature of GPR compared to a generic linear regression of type of Eq. 8 lies in the use of a covariance function for the basis, which gives Eqs. 1-2 the meaning of the estimates, respectively, of the mean and variance of a Gaussian distribution of values of $f(x)$.

The defining equation of GPR, Eq. 1, is thus based on a *square* linear problem with as many training points as basis functions. With $N = M$, Eq. 9 is solvable exactly barring numeric instability of the inversion, and *in this case there is no information inherent in Eq. 9 to optimize the basis*. One then has to resort to heuristics such as MLE for hyperparameter optimization. When $N > M$, however, i.e. with a rectangular version of Eq. 9, it in general cannot be solved exactly, and the residual of a least-squared solution can be used to guide hyperparameter optimization to improve basis completeness (Chan et al., 2012; Manzhos et al., 2011b, 2011a),

$$\min_{\lambda}\left\{N^{-\frac{1}{2}}\left((\boldsymbol{f}' - \boldsymbol{Bc})^T(\boldsymbol{f}' - \boldsymbol{Bc})\right)^{\frac{1}{2}}\right\} \equiv \min_{\lambda}\{rmse_{res}\}$$

(10)

where $\boldsymbol{c} = \boldsymbol{B}^+\boldsymbol{f}'$ and $\boldsymbol{B}^+$ is a Moore-Penrose pseudoinverse (Penrose, 1955) of $\boldsymbol{B}$. In this case one selects (randomly in our case) $M$ from the available $N$ datapoints as basis centers. If one wishes to use $\delta$ in this scheme, it can be introduced into the diagonal of matrix $\Sigma$ of the singular value decomposition of $\boldsymbol{B}$:

$$\boldsymbol{B} = \boldsymbol{U\Sigma V}^T$$
$$\boldsymbol{B}^+ = \boldsymbol{V\Sigma}^+\boldsymbol{U}^T$$

(11)

where $\Sigma^+$ is formed from $\Sigma$ by taking the reciprocal of all the non-zero elements and then transposing (Golub and Van Loan, 1996). However, as we are solving a rectangular problem in the least-squares sense, the parameter $\delta$ is not needed. We found no advantages of using the parameter $\delta$ with a rectangular Eq. 9; the rectangular approach achieves as good a global error without $\delta$ as the traditional GPR with optimal $\delta$ in (Manzhos and Ihara, 2022). As only $M$ points serve as basis centers, we can still use Eq. 2 to approximate the variance of the predicted values with the $\boldsymbol{K}^*$ and $\boldsymbol{K}$ of Eq. 2 computed based on those $M$ points. If $\delta$ is not used, one can use $\boldsymbol{K}^+$ instead of $\boldsymbol{K}^{-1}$ in Eq. 2.

In (Manzhos and Ihara, 2022), we studied the optimization of hyperparameters of GPR from the perspective of basis optimization. In that work, to optimize the parameters of the basis in Eq. 8, we generated with the help of RS-HDMR a reference function $f_{ref}(\boldsymbol{x})$ which was sufficiently close to the original function $f(\boldsymbol{x})$ and such that (i) it avoided overfitting (was well defined with few data) and (ii) a basis which was optimal for $f_{ref}(\boldsymbol{x})$ was also "good enough" for $f(\boldsymbol{x})$. The role of $f_{ref}(\boldsymbol{x})$ which was computable everywhere in space was to enable the calculation of a large test set of points which was then used to optimize basis hyperparameters. In (Manzhos and Ihara, 2022), that approach was also compared to MLE. The results of (Manzhos and Ihara, 2022) are therefore a natural comparison point for the present approach of rectangularization of GPR. We use the same dataset of ab initio samples of the 15-dimensional potential energy surface of UF$_6$ molecule for $f(\boldsymbol{x})$. The data set was made publicly available in the Supporting Information of (Boussaidi et al., 2020), where the plots describing the shape of the function are also given. The PES was sampled in the space of 15 normal mode coordinates with about 55,000 ab initio calculations. The sampling points were distributed in space with the help of Sobol sequence (Sobol', 1967), the point distribution was also made overweigh the data density in the low-energy region, which is an advantageous sampling scheme for non-reactive PESs (Garashchuk and Light, 2001; Manzhos and Carrington, 2016). The values of potential energy range 0-6,629 cm$^{-1}$. Further details of the calculations are given in (Manzhos et al., 2015a); they are not important for the purpose of the present work.

The calculations were performed in Matlab 2021a implementing Eqs. 9-11. A squared exponential kernel of Eq. 6 is used. As all features are normalized to unit standard deviation, it is sufficient to use an isotropic kernel (with a single hyperparameter $l$). Out of the full set of $N_{full}$ = 54,991 available points, we imagine that only a subset of $N$ are available for training the GPR model, out of which we select $M$ as basis centers, to form the matrix $\boldsymbol{B}$ of Eq. 9. We optimize the hyperparameters to minimize $rmse_{res}$ of Eq. 10 and also monitor the rmse on $N_{full} - N$ test points which would have been unavailable in an application. We test whether minimizing $rmse_{res}$ on $N$ points also minimizes the test set rmse. We use test sets which are larger than training sets to monitor global quality of the fit and therefore the quality of the basis and to confirm that minimizing $rmse_{res}$ (which itself becomes possible due to rectangularization) on $N$ points does improve the global quality of GPR.

## 3   Results

The results are summarized in Table 1, and examples of correlation plots between the exact and rectangular GPR-predicted values of the potential energy values for different combinations of $N, M$, and $l$ are shown in Figure 1. Note that all Pearson correlation coefficients are very close to 1 for train and test points. The spread of values visible in the figure optically overstates the extent of error with few points deviating from the diagonal while most points that stay very close to the diagonal visually overlap. In (Manzhos and Ihara, 2022), we manually scanned $l$ and $\delta$ for different $M$ to determine the optimal hyperparameters of the traditional (square) GPR and achievable test set error on the same dataset as used here. The best test set error with the traditional GPR was 36.9 cm$^{-1}$ with 5,000 training points / basis functions and 25.5 cm$^{-1}$ with 10,000 training points / basis functions. The corresponding train set errors were 10.1 cm$^{-1}$ and 13.7 cm$^{-1}$. With rectangular GPR, using 5,000 basis functions and $N$ = 10,000, by minimizing a rectangular residual we obtain an error of about 13 cm$^{-1}$ and a test set error of about 27 cm$^{-1}$ with $l$ = 2.0-2.5. Using 10,000 basis functions and $N$ = 20,000, we obtain a rectangular residual based error of about 13 cm$^{-1}$ and a test set error of about 23 cm$^{-1}$ also with $l$ = 2.0-2.5. With the same number of known points $N$ of rectangular GPR as train points of square GPR, we obtain a similar *global* error as was obtained with square GPR in (Manzhos and Ihara, 2022) *where the hyperparameters were tuned to a large (much larger than N) test set.* Based only on the knowledge of the training point set of 10,000 points and with the help of RS-HDMR, somewhat suboptimal basis hyperparameters could be derived and a global error of 35.1 cm$^{-1}$ could be obtained with square GPR (Manzhos and Ihara, 2022). *With rectangular GPR, the hyperparameters optimal for the global error were obtained based on information contained in N points only.*

Table 1. Root mean square errors obtained with the rectangular Eq. 9 (*rmse$_{res}$*) and on the test set of points, in cm$^{-1}$. *N* is the number of rows in matrix ***B*** and *M* is the number of columns (number of basis functions). Results with the optimal length based on *rmse$_{res}$* for each *N* are highlighted in bold.

| *N* | *M* | *l* | *rmse$_{res}$* | Test rmse |
|---|---|---|---|---|
| 5,000 | 1,000 | 2.0 / 2.5 / 3.0 / 4.0 | 68.4 / 50.8 / 48.7 / 50.3 | 89.1 / 67.3 / 63.9 / 64.1 |
| | 2,000 | 1.5 / 2.0 / 2.5 / 3.5 | 56.6 / 34.5 / 32.1 / 43.4 | 90.0 / 59.0 / 56.3 / 60.0 |
| | 3,000 | 1.5 / **2.0** / **2.5** / 3.0 | 28.3 / **18.6** / **18.4** / 20.5 | 64.5 / **49.8** / **50.8** / 49.5 |
| 10,000 | 2,000 | 1.5 / 2.0 / 2.5 / 3.0 | 59.7 / 41.9 / 37.6 / 40.6 | 77.6 / 53.4 / 50.1 / 52.6 |
| | 3,000 | 1.5 / 2.0 / 2.5 / 3.0 | 36.2 / 27.0 / 26.7 / 29.2 | 52.0 / 40.5 / 40.6 / 42.9 |
| | 5,000 | 1.5 / **2.0** / **2.5** / 3.0 | 19.4 / **13.2** / **13.5** / 15.7 | 37.6 / **27.2** / **26.9** / 28.3 |
| | 7,000 | 1.5 / 2.0 / 2.5 / 3.0 | 14.8 / 10.8 / 11.9 / 14.7 | 40.1 / 31.5 / 28.2 / 27.1 |
| 15,000 | 3,000 | 1.5 / 2.0 / 2.5 / 3.0 | 39.4 / 29.7 / 28.2 / 32.2 | 50.1 / 39.2 / 37.5 / 41.2 |
| | 5,000 | 1.5 / 2.0 / 2.5 / 3.0 | 22.1 / 15.3 / 15.3 / 20.6 | 32.5 / 23.9 / 23.2 / 29.6 |
| | 7,500 | 1.5 / **2.0** / **2.5** / 3.0 | 17.9 / **12.9** / **13.3** / 17.3 | 33.3 / **25.4** / **24.2** / 24.6 |
| | 10,000 | 1.5 / 2.0 / 2.5 / 3.0 | 14.4 / 10.5 / 12.3 / 16.8 | 36.7 / 28.9 / 24.6 / 24.3 |
| 20,000 | 3,000 | 1.5 / 2.0 / 2.5 / 3.0 | 40.4 / 29.9 / 28.8 / 33.6 | 49.5 / 36.7 / 35.7 / 39.8 |
| | 5,000 | 1.5 / 2.0 / 2.5 / 3.0 | 22.9 / 16.5 / 16.4 / 22.2 | 31.0 / 22.6 / 22.1 / 28.4 |
| | 10,000 | 1.5 / **2.0** / **2.5** / 3.0 | 14.0 / **12.2** / **13.6** / 17.7 | 31.2 / **24.6** / **22.4** / 23.0 |
| | 15,000 | 1.5 / 2.0 / 2.5 / 3.0 | 11.0 / 9.3 / 12.3 / 18.1 | 39.9 / 31.0 / 22.9 / 23.4 |

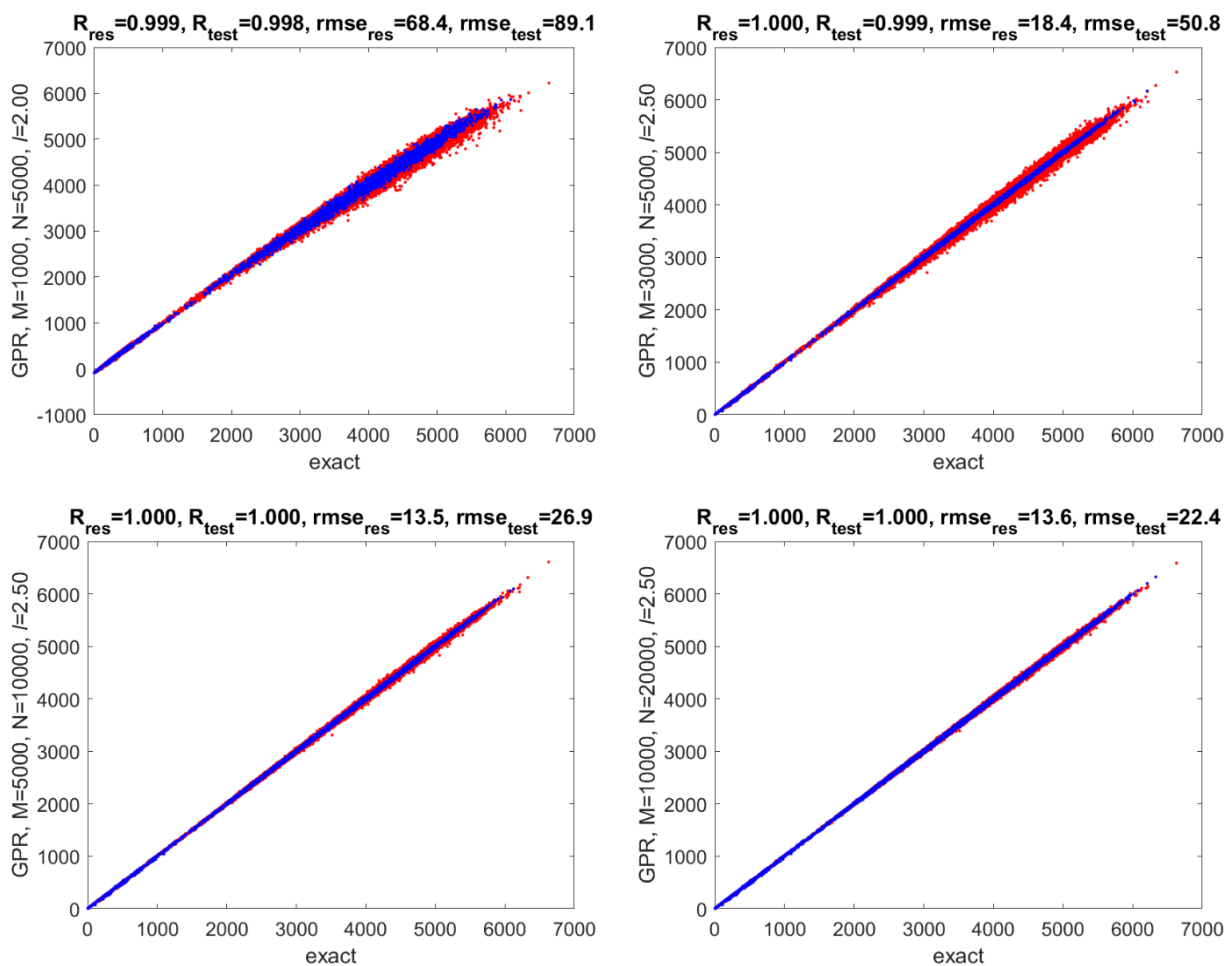

Figure 1. Example correlation plots between the exact and GPR-predicted values of the potential energy, in cm$^{-1}$, for different *M, N, l* combinations. Blue points are for the training set and red points for the test set. The values of Pearson correlation coefficient for the *N* points used in the training and the $N_{all}$-*N* test points are shown on the plots as well as corresponding rmse values.

Important for this is the observation that choosing *l* by minimizing the residual of the rectangular matrix equation also minimizes the global error proxied here by the error on a large test set (much larger than the training set). The ratio of *N* to *M* that provides a good predictive power of optimal *l* is roughly 2. When *M* is more than about 1/2 of *N*, there is a danger of overfitting and the trend in *rmse*$_{res}$ as a function of *l* may no longer reflect the trend in test rmse. To put this in perspective, a ratio of about 2-3 of the number of points to the number of basis functions was found to be optimal for the rectangular collocation approach (to build a basis representation of vibrational wavefunctions) with Gaussian and other bell-like basis functions (Kamath and Manzhos, 2018; Manzhos et al., 2018). One can therefore use much fewer basis functions than known points, which decreases the cost of fitting and calling the model. Note that optimal *l* values with the "square" GPR (which were about 3.0-3.5) were larger than those we obtain with the rectangular GPR, and optimal $\delta$ was on the order of $1 \times 10^{-4\ldots-5}$. These values of $\delta$ were close to the threshold of numeric stability (calculations with $\delta$ values below $1 \times 10^{-6}$ were unstable) (Manzhos and Ihara, 2022). With rectangular GPR, we dispense with $\delta$ altogether which simplifies hyperparameter optimization (our tests with $\delta$ added to the diagonal of $\Sigma$ showed an effect on rmse but no additional improvements in the test set rmse).

In Figure 2, we show the three-sigma confidence intervals for selected cases (*M, N, l* combinations). To make the figure readable, the plots are done for random subsets of 1,000 test points of which 500 show confidence intervals, to visually appreciate the relation of computed confidence intervals to the actual spread of the data. We observe that the confidence intervals depend strongly on *l* and may or may not correspond to the actual error of prediction. This behavior is the same as in the traditional square GPR (Boussaidi et al., 2020). We caution that confidence intervals computed with Eq. 2 should not automatically be used as error bars. See notes to this effect in (Boussaidi et al., 2020; Ren et al., 2021).

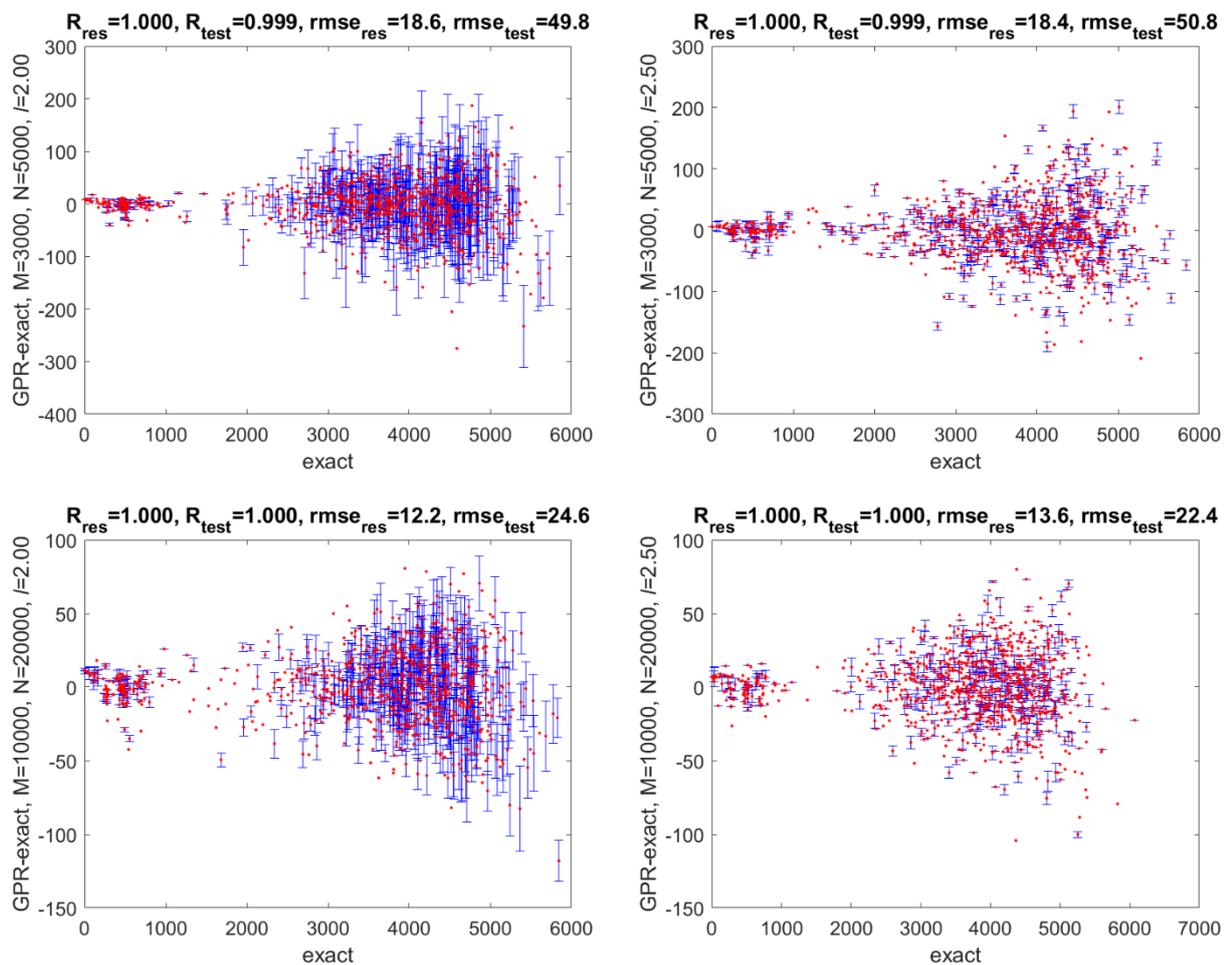

Figure 2. Difference between the GPR-computed and exact values of the potential energy, in cm$^{-1}$, for selected *M, N, l* combinations. For readability, a random subset of 1,000 test points is shown of which 500 show three-sigma confidence intervals, to visually appreciate the relation of confidence intervals computed with Eq. 2 to the actual spread of the data.

# 4 Conclusions

We explored rectangularization of the defining equation of Gaussian process regression to address the important problem of optimization of hyperparameters which determines success or failure of the application of the GPR method. Such optimization is difficult when data are scarce (data density is low), with failures of the MLE method likely and documented. Recognizing that GPR is equivalent to a linear regression with a basis formed by covariance functions, we proposed a hyperparameter optimization approach based on the quality of the least-squares solution of a rectangular matrix equation of size $N \times M$, where $M$ is the basis size and $N > M$ is the number of points used to tune hyperparameters. We demonstrated the advantages of the approach on the example of fitting a 15-dimensional potential energy surface of $UF_6$ molecule and showed that it allows effective hyperparameter tuning even with very scarce data. With less than 10,000 data points the method was able to determine hyperparameters which were optimal for the global quality of the function (proxied by a large test set). The rectangular GPR makes unnecessary the noise parameter typically added to the diagonal of the covariance matrix to improve stability and generalization, which simplifies hyperparameter optimization. One can use much fewer basis functions than known points, which decreases the cost of fitting and calling the model. Using about a half of known points as basis centers can be recommended, and good results are obtainable when only a third of the data points are used for that purpose.